\documentclass[11pt]{article}
\usepackage{mathrsfs}
\usepackage{amsmath,amsfonts,amssymb,rotating,amsthm}
\usepackage{hyperref}
\usepackage{array}
\usepackage{breqn}
\allowdisplaybreaks
\def\qed{\nopagebreak\hfill{\rule{4pt}{7pt}}}
\def\proof{\noindent {\it{Proof.} \hskip 2pt}}
\newtheorem{theo}{Theorem}[section]

\newtheorem{lemm}[theo]{Lemma}

\newtheorem{conj}[theo]{Conjecture}
\theoremstyle{remark}

\numberwithin{equation}{section}

\newdimen\Squaresize \Squaresize=11pt
\newdimen\Thickness \Thickness=0.7pt
\def\Square#1{\hbox{\vrule width \Thickness
   \vbox to \Squaresize{\hrule height \Thickness\vss
    \hbox to \Squaresize{\hss#1\hss}
   \vss\hrule height\Thickness}
\unskip\vrule width \Thickness} \kern-\Thickness}

\def\Vsquare#1{\vbox{\Square{$#1$}}\kern-\Thickness}

\def\moins{\raise 1pt\hbox{{$\scriptstyle -$}}}

\begin{document}

\begin{center}
{\large \bf  Higher Order Tur\'{a}n Inequalities

for the Partition Function}
\end{center}

\begin{center}
William Y.C. Chen$^{1}$, Dennis X.Q. Jia$^{2}$ and Larry X.W. Wang$^{3}$\\[8pt]
Center for Combinatorics, LPMC\\
Nankai University\\
 Tianjin 300071, P. R. China\\[6pt]
Email: $^{1}${\tt
chen@nankai.edu.cn, $^{2}$jiaxiaoqiang2017@gmail.com, $^{3}$wsw82@nankai.edu.cn\tt }
\end{center}

\vspace{0.3cm} \noindent{\bf Abstract.}

The Tur\'{a}n inequalities and the higher order Tur\'{a}n inequalities arise in the study of  Maclaurin coefficients of an entire function in the Laguerre-P\'{o}lya class. A real sequence $\{a_{n}\}$  is said to satisfy the Tur\'{a}n inequalities if for $n\geq 1$, $a_n^2-a_{n-1}a_{n+1}\geq 0$. It  is said to satisfy the higher
 order Tur\'{a}n inequalities if for $n\geq 1$,
$4(a_{n}^2-a_{n-1}a_{n+1})(a_{n+1}^2-a_{n}a_{n+2})-(a_{n}a_{n+1}-a_{n-1}a_{n+2})^2\geq 0$.
A sequence satisfying the Tur\'an inequalities is also called
 log-concave. For the partition function $p(n)$,  DeSalvo and Pak
showed that for $n>25$, the sequence $\{ p(n)\}_{n> 25}$ is
log-concave, that is, $p(n)^2-p(n-1)p(n+1)>0$ for $n> 25$.
 It was conjectured by Chen  that  $p(n)$
 satisfies the higher order Tur\'{a}n inequalities for $n\geq 95$.
In this paper, we prove this conjecture by using the Hardy-Ramanujan-Rademacher formula to
derive an upper bound and a lower bound for $p(n+1)p(n-1)/p(n)^2$.
Consequently, for $n\geq 95$, the Jensen polynomials $g_{3,n-1}(x)=p(n-1)+3p(n)x+3p(n+1)x^2+p(n+2)x^3$ have only real zeros.
We  conjecture that for any positive integer
$m\geq 4$ there exists an integer $N(m)$ such that for $n\geq N(m) $, the polynomials $\sum_{k=0}^m {m\choose k}p(n+k)x^k$ have only real zeros.
This conjecture was independently posed by Ono.

\\

\\

\section{Introduction}

The objective of this paper is to prove the higher order Tur\'{a}n inequalities for the partition function $p(n)$ when $n \geq 95$, as  conjectured   in \cite{Chen3}.
The Tur\'{a}n inequalities and the higher order Tur\'an inequalities
are related to  the Laguerre-P\'{o}lya class of real entire functions (cf. \cite{Dimitrov, Szego}).
A sequence $\{a_n\}$ of real numbers is said to satisfy the Tur\'an inequalities
if
\begin{equation}\label{T}
   a_n^2-a_{n-1}a_{n+1}\geq 0,
\end{equation}
for $n\geq 1$. The inequalities   \eqref{T} are also called the Newton inequalities (cf. \cite{Craven, Csordas, Niculescu, Wagner}). We say that a sequence $\{a_n\}$ satisfies the higher
 order Tur\'{a}n inequalities or cubic Newton inequalities (cf. \cite{Niculescu}) if for $n\geq 1$,
 \begin{equation}\label{HT}
   4(a_{n}^2-a_{n-1}a_{n+1})(a_{n+1}^2-a_{n}a_{n+2})-(a_{n}a_{n+1}-a_{n-1}a_{n+2})^2\geq 0.
\end{equation}
A real entire function
\begin{eqnarray} \label{1}
 \psi(x)=\sum_{k=0}^{\infty} \gamma_k\frac{x^k}{k!}
 \end{eqnarray} is said to be in the Laguerre-P\'{o}lya class, denoted
  $\psi(x) \in \mathcal{LP}$, if it can be represented in the form
\begin{eqnarray*}
\psi(x)=c x^m e^{-\alpha x^{2}+\beta x} \prod_{k=1}^{\infty}\left(1+x/x_{k}\right)e^{- x/x_k},
\end{eqnarray*}
where $c$, $\beta$, $x_{k}$ are real numbers, $\alpha \geq 0$, $m$ is a nonnegative integer and $\sum x^{-2}_{k}< \infty$. These functions are the only ones which are uniform limits of polynomials whose zeros are real. We refer to \cite{Levin} and \cite{Rahman} for the background on the theory of the $\mathcal{LP}$ class.

Jensen \cite{Jensen} proved
 that
a real entire function $\psi(x)$ belongs to $\mathcal{LP}$ class if and only if for any positive integer $n$, the $n$-th associated Jensen polynomial
\begin{equation}
g_n(x)=\sum_{k=0}^{n}{ n\choose{k}} \gamma_k x^k
\end{equation}
has only real zeros.
More properties
 of the Jensen polynomials can be found in  \cite{Craven, Csordas, Csordas2}.

P\'olya and Schur \cite{Polya} also obtained the above result based on multiplier sequences of the second kind.
A real sequence $\{\gamma_k\}_{k\geq 0} $  is called a multiplier sequence of the second kind if for any nonnegative integer $n$ and every real polynomial $\sum_{k=0}^{n} a_k x^k$ with only real zeros  of the same sign,
the polynomial $\sum_{k=0}^{n} a_k \gamma_k x^k$ has only real zeros.   P\'olya and Schur \cite{Polya}  proved that each multiplier sequence of the second kind  satisfies the Tur\'{a}n inequalities. Moreover, they showed that a real entire function $\psi(x)$  belongs to the $\mathcal{LP}$ class if and only if its Maclaurin coefficient sequence is a multiplier sequence of the second kind. This implies that the Maclaurin coefficients of $\psi(x)$ in the $\mathcal{LP}$ class satisfy the Tur\'{a}n inequalities
\begin{eqnarray} \label{1.7}
\gamma_{k}^2-\gamma_{k-1}\gamma_{k+1}\geq 0
\end{eqnarray}
for $k\geq 1$.
In fact, \eqref{1.7} is a consequence of
another property of the $\mathcal{LP}$ class due to P\'olya and Schur \cite{Polya}:
Let $\psi(x)$ be a real entire function in the $\mathcal{LP}$ class.  For any nonnegative integer $m$,
 the $m$-th derivative $\psi^{(m)}$ of $\psi(x)$ also belongs to
the $\mathcal{LP}$ class.
It is readily seen that the $n$-th Jensen polynomial associated with $\psi^{(m)}$ is
\begin{eqnarray}\label{9}
g_{n,m}(x)=\sum_{k=0}^{n}{n\choose k} \gamma_{k+m} x^k,
\end{eqnarray}
and hence it
 has only real zeros for any nonnegative integers $n$ and $m$.
In particular, taking $n=2$, for any nonnegative integer $m$,
 the real-rootedness of $g_{2,m}(x)$
 implies the inequality \eqref{1.7}, that is, the discriminant $4(\gamma_{m+1}^2-\gamma_{m}\gamma_{m+2})$ is nonnegative.

Dimitrov \cite{Dimitrov} observed that for  a real entire function $\psi(x)$   in the $\mathcal{LP}$ class,
the Maclaurin coefficients satisfy the higher
order Tur\'{a}n inequalities
\begin{equation}\label{8}
4(\gamma_{k}^2-\gamma_{k-1}\gamma_{k+1})(\gamma_{k+1}^2-\gamma_{k}
\gamma_{k+2})-(\gamma_{k}\gamma_{k+1}-\gamma_{k-1}\gamma_{k+2})^2\geq0
\end{equation}
 for $k\geq 1$.
This fact follows from a theorem of
Ma\v{r}\'{\i}k \cite{Marik} stating that if a real polynomial
   \begin{equation}
 \sum_{k=0}^{n} {n\choose k} a_k x^k
   \end{equation}
   of degree $n\geq 3$ has only real zeros, then
 $a_0, a_1, \ldots, a_n$  satisfy the higher
 order Tur\'{a}n inequalities.

As noted in \cite{Chen3}, for
 $k=1$, the polynomial
in \eqref{8} coincides with an invariant
\begin{equation}\label{equ-dis-cub}
\nonumber I(a_0,a_1,a_2,a_3)=3a_1^2 a_2^2-4a_1^3 a_3-4a_0 a_2^3-a_0^2 a_3^2+6a_0a_1 a_2a_3
\end{equation}
 of the cubic binary form
 \begin{eqnarray}\label{equ-cub-form}
a_3 x^3+3a_2 x^2 y+3 a_1 x y^2+a_0 y^3.
\end{eqnarray}
In other words, the above invariant $I(a_0,a_1,a_2,a_3)$ can be rewritten as
\begin{equation}
 I(a_0,a_1,a_2,a_3)=  4(a_1^2 -a_0 a_2)(a_2^2- a_1 a_3)-(a_1 a_2- a_0 a_3)^2.
\end{equation}
We refer to
Hilbert \cite{Hilbert}, Kung and Rota \cite{KR} and Sturmfels \cite{Sturmfels} for the background on the invariant theory of
 binary forms.
Notice that   $27I(a_0, a_1, a_2, a_3)$ is  the discriminant of the cubic polynomial $a_3 x^3+3a_2 x^2 +3 a_1 x +a_0 $.
A  cubic polynomial with real coefficients
 has only real zeros if and
only if its discriminant is nonnegative \cite[p.42]{Cajori}.

Recall that for a real entire function $\psi(x)$ in the $\mathcal{LP}$ class,
its $m$-th derivative  $\psi^{(m)}$ is also a real entire function  in the  $\mathcal{LP}$ class.
Thus the real-rootedness of  the  cubic Jensen polynomial  $g_{3,m}(x)$ associated with $\psi^{(m)}$  implies the higher order Tur\'an inequalities \eqref{8} of Dimitrov, that is, the discriminant $27I(\gamma_{m}, \gamma_{m+1}, \gamma_{m+2}, \gamma_{m+3})$ is nonnegative.

Real entire functions   in the $\mathcal{LP}$ class with
nonnegative Maclaurin coefficients also received much attention.
Aissen, Schoenberg and Whitney \cite{Aissen} proved that if $\psi(x)$
 is a real entire function in the $\mathcal{LP}$ class with
nonnegative Maclaurin  coefficients, then the
 the  sequence
 $\{\gamma_k/k!\}$ associated with $\psi(x)$ forms a {P\'olya frequency sequence}.
 An infinite sequence $\{a_n\}_{n\geq0}$ of nonnegative numbers is called a {P\'olya frequency sequence} (or a $PF$-sequence) if the matrix $(a_{i-j})_{i,j\geq0}$ is a totally positive matrix, where $a_n=0$ if $n<0$, that is, all minors of $(a_{i-j})_{i,j\geq0}$  have nonnegative determinants.
For more properties of totally positive matrices  and $PF$-sequences can be found in \cite{Craven2, Karlin}.

The $\mathcal{LP}$ class is closely related to the Riemann hypothesis.
Let $\zeta$ denote the Riemann zeta-fucntion and $\Gamma$ be the gamma-function.
The Riemman $\xi$-function is defined by
\begin{eqnarray}
\xi(iz)=\frac{1}{2}\left(z^2-\frac{1}{4}\right)
\pi^{-z/2-1/4}\Gamma\left(\frac{z}{2}+\frac{1}{4}\right)\zeta\left(z+\frac{1}{2}\right),
\end{eqnarray}
see, for example, Boas \cite[p.24]{Boas}.
It is well known that the Riemann hypothesis holds  if and only if
the  Riemann $\xi$-function belongs to the  $\mathcal{LP}$ class.
Hence, if the Riemann hypothesis is true, then the Maclaurin coefficients of the Riemann $\xi$-function satisfy  both the Tur\'{a}n inequalities
and the higher order Tur\'{a}n inequalities.
Csordas, Norfolk and Varga \cite{Csordas}  proved that the coefficients of the Riemann $\xi$-function satisfy the Tur\'{a}n inequalities, confirming a conjecture of P\'{o}lya. Dimitrov and Lucas \cite{Dimitrov2} showed that the coefficients of the Riemann $\xi$-function satisfy the higher order Tur\'{a}n inequalities without resorting to the Riemann hypothesis.

Let us now turn to the partition function. A partition of a  positive integer $n$ is a nonincreasing sequence $(\lambda_1,\lambda_2,\ldots, \lambda_r)$ of positive integers such that $ \lambda_1+\lambda_2+\cdots+\lambda_r=n$. Let $p(n)$ denote the number of partitions of $n$.
A sequence $\{a_k\}_{k\geq 0}$ satisfying the Tur\'{a}n inequalities, that is, $a_k^2-a_{k-1}a_{k+1}\geq 0$ for $k\geq 1$, is also called {log-concave}. DeSalvo and Pak \cite{Desalvo} proved the log-concavity of the partition function
for $n> 25$ as well as  the following  inequality as conjectured
in \cite{Chen}:
\begin{eqnarray}
 \frac{p(n-1)}{p(n)}\left(1+\frac{1}{n}\right)>\frac{p(n)}{p(n+1)}
\end{eqnarray}
for $n\geq 2$.
DaSalvo and Pak  also conjectured that for $n\geq 45$,
\begin{eqnarray}\label{DP}	 \frac{p(n-1)}{p(n)}\left(1+\frac{\pi}{\sqrt{24}n^{3/2}}\right)>\frac{p(n)}{p(n+1)}.
\end{eqnarray}
Chen, Wang and Xie \cite{Chen2} gave an affirmative answer to this conjecture.

It was conjectured  in \cite{Chen3} that for large $n$,
 the partition function $p(n)$ satisfies many inequalities pertaining to
 invariants
of a binary form. Here we mention two of them.

\begin{conj}\label{conj1}
For $n\geq 95$, the higher order Tur\'an inequalities
 \begin{eqnarray}\label{1.11}
4(a_{n}^2-a_{n-1}a_{n+1})(a_{n+1}^2-a_{n}
a_{n+2})-(a_{n}a_{n+1}-a_{n-1}a_{n+2})^2\geq0
\end{eqnarray}
 hold for $a_n=p(n)$.
\end{conj}

The following conjecture is a higher order analogue of
\eqref{DP}.

\begin{conj}
Let
\begin{eqnarray}\label{1.8}
u_n=\frac{p(n+1)p(n-1)}{p(n)^2}.
\end{eqnarray}
For $n\geq 2$,
	\begin{eqnarray*}
4\left(1-u_n\right)\left(1-u_{n+1}\right)
<\left(1+\frac{\pi}{\sqrt{24}n^{3/2}}\right)\left(1-u_n u_{n+1}\right)^2.
	\end{eqnarray*}
\end{conj}

The objective  of this paper is to prove Conjecture \ref{conj1}.
In fact, we shall prove the following equivalent form.

\begin{theo} \label{Theorem 2.2}
	Let $u_n$ be defined as in \eqref{1.8}. For $n \geq 95$,
	\begin{eqnarray}\label{an}
		4(1-u_n)(1-u_{n+1})-(1-u_n u_{n+1})^2>0.
	\end{eqnarray}
\end{theo}

The above theorem can be restated as follows.

\begin{theo}
For $n\geq 95$, the cubic polynomial
\[
p(n-1)+3p(n)x+3p(n+1)x^2+p(n+2)x^3
\]
has three distinct real zeros.
\end{theo}

In general, we propose
 the following conjecture.

\begin{conj}\label{conj3}
For any positive integer $m\geq 4$, there exists a positive  integer $N(m)$ such that for any $n\geq N(m)$, the polynomial
\[\sum_{k=0}^m {m\choose k}p(n+k)x^k\]
 has only real zeros.
 \end{conj}

 The above conjecture was independently proposed by Ono \cite{Ono}.
 It was recently announced by Ono that he and Zagier have
 proved this conjecture.

\section{Bounding $u_n$}

In this section, we give an upper  bound and a lower bound for
\[ u_n=\frac{p(n+1)p(n-1)}{p(n)^2},\]
as defined in \eqref{1.8}.
DeSalvo and Pak \cite{Desalvo}   proved that for $n>25$,
\[1-\frac{1}{n+1}<u_n<1 .\]
On the other hand, Chen, Wang and Xie \cite{Chen2}  showed that for $n\geq45$,
\[1-\frac{\pi}{\sqrt{24}n^{3/2}+\pi}<u_n.\]
Nevertheless, the above  bounds for $u_n$ are not sharp enough for the purpose of
proving Theorem \ref{Theorem 2.2}.
To state our bounds for $u_n$, we adopt the following notation as used in \cite{Lehmer2}:
\begin{equation}\label{mu}
    \mu(n)=\frac{\pi}{6}   \sqrt{24 n-1}.
\end{equation}
For convenience, let
\begin{eqnarray}\label{xyzw}
	x=\mu(n-1),\ y=\mu(n),\ z=\mu(n+1),\ w=\mu(n+2).
\end{eqnarray}
Define
\begin{eqnarray}\label{f}
f(n)=e^{x-2 y+z} \frac{\left(x^{10}-x^9-1\right) y^{24} \left(z^{10}-z^9-1\right)}{x^{12} \left(y^{10}-y^9+1\right)^2 z^{12}},
\end{eqnarray}
\begin{eqnarray}\label{g}
g(n)= e^{x-2 y+z} \frac{\left(x^{10}-x^9+1\right) y^{24} \left(z^{10}-z^9+1\right)}{x^{12} \left(y^{10}-y^9-1\right)^2 z^{12}}.
\end{eqnarray}

Then we have the following bounds for $u_n$.

\begin{theo} \label{Lemmax1}
For $n\geq 1207$,
	\begin{eqnarray}
	&&f(n)<u_n<g(n). \label{2.1}
	\end{eqnarray}
\end{theo}

In order to give a proof of Theorem \ref{Lemmax1},
we need the following  upper bound and   lower bound for $p(n)$.

\begin{lemm} \label{Lemma2-1}
Let
\begin{align*}
    &B_1(n)=\frac{\sqrt{12}e^{\mu(n)}}{24n-1}\left(1-\frac{1}{\mu(n)}
-\frac{1}{\mu(n)^{10}}\right),\\[9pt]
&B_2(n)=\frac{\sqrt{12}e^{\mu(n)}}{24n-1}\left(1-\frac{1}{\mu(n)}
+\frac{1}{\mu(n)^{10}}\right),
\end{align*}
then for $n\geq 1206$,
\begin{equation}\label{bound-p}
B_1(n)<p(n)<B_2(n).
\end{equation}
\end{lemm}

The proof of Lemma \ref{Lemma2-1} relies on the
 Hardy-Ramanujan-Rademacher formula \cite{Hardy, Rademacher}  for
 $p(n)$ as well as Lehmer's error bound for the remainder of this formula.
The Hardy-Ramanujan-Rademacher formula reads
\begin{align}\label{1.5}
\nonumber p(n)=&\frac{\sqrt{12}}{24n-1}\sum^{N}_{k=1}\frac{A_{k}(n)}{\sqrt{k}}
\left[\left(1-\frac{k}{\mu(n)}\right)e^{\mu(n)/k}
+\left(1+\frac{k}{\mu(n)}\right)e^{-\mu(n)/k}\right]\\[9pt]
&\qquad +R_{2}(n,N),
\end{align}
where $A_{k}(n)$ is an arithmetic function and  $R_{2}(n,N)$ is the remainder term,
see, for example,  Rademacher \cite{Rademacher}. Lehmer \cite{Lehmer1,Lehmer2} gave the following error bound:
\begin{equation}\label{R2}
|R_{2}(n,N)|<\frac{\pi^{2}N^{-2/3}}{\sqrt{3}}\left[\left(\frac{N}{\mu(n)}\right)^{3}
\sinh \frac{\mu(n)}{N}+\frac{1}{6}-\left(\frac{N}{\mu(n)}\right)^{2}\right],
\end{equation}
which is valid for all positive integers $n$ and $N$.

\noindent
{\it Proof of Lemma \ref{Lemma2-1}. }
Consider  the Hardy-Ramanujan-Rademacher formula   \eqref{1.5} for  $N=2$, and note that $A_1(n)=1$ and $A_2(n)=(-1)^n$  for any positive integer $n$. Hence \eqref{1.5} can be rewritten as
\begin{align}\label{pn1}
	{p(n)} ={\frac{\sqrt{12}e^{\mu(n)}}{24n-1}}\left( 1-\frac{1}{\mu(n)}+T(n)\right),
\end{align}
where
	\begin{align}
T(n)=&\frac{(-1)^n}{\sqrt{2}} \left(\left(1-\frac{2}{\mu(n)}\right)e^{-\mu(n)/2}+
\left(1+\frac{2}{\mu(n)}\right)e^{-3\mu(n)/2}\right)\nonumber\\[6pt]
& \qquad+\left(1+\frac{1}{\mu(n)}\right) e^{-2\mu(n)} +(24n-1)R_2(n,2)/\sqrt{12}e^{\mu(n)}.\label{Tn1}
\end{align}
In order to prove \eqref{bound-p},
we proceed to use Lehmer's error bound in \eqref{R2}  to show that for $n>1520$,
\begin{align}\label{tn1}
|T(n)|<\frac{1}{\mu(n)^{10}}.
\end{align}
Let us begin with the first three terms in   \eqref{Tn1}.
Evidently, for $n\geq 1$,
\[0<\frac{1}{\mu (n)}<\frac{1}{2},\]
so that
\begin{align}
&\left(1-\frac{2}{\mu (n)}\right) e^{{-\mu (n)}/{2} }<e^{{-\mu(n)}/{2}},\label{Y1}\\[9pt]
&\left(1+\frac{2}{\mu  (n)}\right) e^{{-3\mu(n)}/{2}}<2e^{{-3\mu(n)}/{2}},\label{Y2}\\[9pt]
&\left(1+\frac{1}{\mu(n)}\right) e^{-2\mu(n)}<2 e^{-2\mu(n)}\label{Y3}.
\end{align}

As for the last term in \eqref{Tn1}, we claim that for $n>350$,
\begin{equation}\label{Y4}
\left|\frac {(24 n -1)R_ 2 (n, 2)} { {\sqrt {12} e^{\mu (n)}}}\right|<e^{ {-\mu (n)}/ {2} }.
\end{equation}
Applying \eqref{R2} with $N=2$, we obtain that
\begin{align}
&\left|\frac {({24 n -1})R_ 2 (n, 2)} { {\sqrt {12} e^{\mu (n)}}}\right| = \left|\frac {36\mu (n)^2 R_ 2 (n, 2)} {\sqrt {12}\pi^2 e^{\mu (n)}}\right|\nonumber\\[9pt]
&<\frac {\mu (n)^2 e^{-\mu (n)}} {2^{2/3}} + \frac {12\sqrt[3] {2} e^{ {-\mu (n)}/ {2} }} {\mu (n)} - \frac {12\sqrt[3] {2} e^{- {3\mu (n)} /{2}}} {\mu (n)} - 12\sqrt[3] {2} e^{-\mu (n)}\nonumber\\[9pt]
&<\frac {\mu (n)^2 e^{-\mu (n)}} {2^{2/3}} + \frac {12\sqrt[3] {2} e^{ {-\mu (n)}/ {2} }} {\mu (n)}\nonumber\\[9pt]
&<  \frac {24 e^{{-\mu (n)}/ {2} }} {\mu (n)}+\mu (n)^2 e^{-\mu (n)} \label{ineq-2-R}.
\end{align}
To bound the first term in \eqref{ineq-2-R},
we find that for $n> 350$,
\begin{align}\label{ineq2-1}
\frac {24 e^{{-\mu (n)}/ {2} }} {\mu (n)}< \frac{e^{ -\mu (n)/2 }}{2},
\end{align}
which simplifies to
\begin{align}\label{ineq2-mu}
\mu(n)=\frac{\pi}{6}   \sqrt{24 n-1}>48,
\end{align}
which is  true for $n> 350$.
Concerning the second term in \eqref{ineq-2-R},
it will be shown that for $n>22$,
\begin{equation}\label{ineq2-2}
    \mu (n)^2 e^{-\mu (n)}<\frac{e^{ -\mu (n)/2 }}{2},
\end{equation}
which can be rewritten as
\begin{equation}\label{ineq2-3}
\frac{e^{\mu(n)/4}}{\mu(n)/4}>4\sqrt{2}.
\end{equation}
Let
\begin{equation}\label{ineq-2-F}
F(t)=\frac{e^{t}}{t}.
\end{equation}
Since $F'(t)={e^{t}(t-1)}/{t^2}>0$ for $t>1$, $F(t)$ is increasing for $t>1$.
Thus,
\[F\left(\frac{\mu(n)}{4}\right)=\frac{e^{\mu(n)/4}}{\mu(n)/4}>F(3)
=\frac{e^3}{3}>4\sqrt{2}.\]
Here we have used the fact that  for $n> 22$, ${\mu(n)}/{4}>3$.
This proves \eqref{ineq2-3}.
Applying the estimates \eqref{ineq2-1} and \eqref{ineq2-2} to \eqref{ineq-2-R}, we
reach \eqref{Y4}.

Taking all the above estimates into account,
 we deduce that for $n> 350$,
\begin{align}\label{T1}
|T(n)|<6e^{-\mu(n)/2}.
\end{align}
To obtain \eqref{tn1}, we have only to  show that for
$n>1520$,
\begin{eqnarray}\label{ineq-mu}
    6e^{-\mu(n)/2}< \frac{1}{\mu(n)^{10}},
\end{eqnarray}
which can be recast as
\begin{equation}\label{ineq-mu2}
 \frac{e^{\mu(n)/20}}{\mu(n)/20}> 20\sqrt[10]{6}.
\end{equation}
Since ${\mu(n)}/{20}>5$ for $n>1520$, 
by  the monotone property of $F(t)$, we have  for $n>1520$,
\[F\left(\frac{\mu(n)}{20}\right)= \frac{e^{\mu(n)/20}}{\mu(n)/20}>F(5)=\frac{e^5}{5}>20 \sqrt[10]{6},\]
as asserted by \eqref{ineq-mu2}.
Thus \eqref{tn1} follows from \eqref{T1} and \eqref{ineq-mu}.
In other words, for  $n>1520$,
\begin{equation}\label{Ineq-2-T_n}
    -\frac{1}{\mu(n)^{10}}<T(n)<\frac{1}{\mu(n)^{10}}.
\end{equation}
Substituting \eqref{pn1} into \eqref{Ineq-2-T_n}, we see that \eqref{bound-p} holds
for $n>1520$. It is routine to check that \eqref{bound-p} is true for $1206\leq n \leq 1520$, and hence the proof is complete.
\qed

We are now ready to prove Theorem \ref{Lemmax1} by Lemma \ref{Lemma2-1}.

\noindent
{\it Proof of Theorem \ref{Lemmax1}. }
Since $B_1(n)$ and $B_2(n)$ are all positive for  $n\geq 1$,
using the bounds for $p(n)$ in \eqref{bound-p}, we find that for $n\geq1207$,
\begin{equation*}
\frac{B_1(n-1)B_1(n+1)}{B_2(n)^2}
<\frac{p(n-1)p(n+1)}{p(n)^2}<\frac{B_2(n-1)B_2(n+1)}{B_1(n)^2}.
\end{equation*}
This proves \eqref{2.1}.
\qed

\section{An inequality on $f(n)$ and $g(n)$}

In this section, we establish an inequality between $f(n)$ and $g(n+1)$
which will be used in the proof of Theorem \ref{Theorem 2.2}.

\begin{theo} \label{Lemmax2}
For $n\geq 2$,
	\begin{align}\label{ineq3-1}
	g(n+1)<f(n)+\frac{110}{\mu(n-1)^5}.
	\end{align}
\end{theo}

\proof
Recall that
\begin{equation}\label{mu-2}
    \mu(n)=\frac{\pi\sqrt{24n-1}}{6},
\end{equation}
and
\begin{equation}\label{x-y-z-w}
    x=\mu(n-1),\ y=\mu(n),\ z=\mu(n+1),\ w=\mu(n+2).
\end{equation}
Let
\begin{align}\label{varphi}
    \alpha(t)=t^{10}-t^{9}+1,\ \ \beta(t)=t^{10}-t^{9}-1.
\end{align}
By the definitions of $f(n)$ and $g(n)$ as given in \eqref{f} and \eqref{g},
we find that
\begin{align}\label{last}
f(n)x^5-g(n+1)x^5+110 =\frac{-e^{w+y-2z}t_1+e^{z+ x-2y}t_2+110t_3}{t_3},
\end{align}
where
\begin{align}
&t_1= x^{12} z^{36} \alpha(y)^3\alpha(w),\label{t-1}\\[9pt]
&t_2= y^{36} w^{12}\beta(x) \beta(z)^3,\label{t-2}\\[9pt]
&t_3= x^7 y^{12} z^{12} w^{12} \alpha(y)^2\beta(z)^2. \label{t-3}
\end{align}
Since $t_3>0$ for $n\geq 2$, \eqref{ineq3-1} is equivalent to
\begin{equation}\label{eq4.1}
    -e^{w+y-2z}t_1+e^{z+ x-2y}t_2+110t_3>0,
\end{equation}
for $n\geq 2$.
To verify \eqref{eq4.1}, we proceed to estimate $t_1, t_2, t_3, e^{w+y-2z}$ and
 $e^{x-2y+z}$ in terms of $x$. Noting that for $n\geq 2$,
\begin{align}\label{eq}
y=\sqrt{x^2+{\frac{2\pi^2}{3}}},\ z=\sqrt{x^2+\frac{4\pi^2}{3}},\ w=\sqrt{x^2+2\pi^2},
\end{align}
we obtain the following expansions
\begin{align*}
    &y=x+\frac{\pi ^2}{3 x}-\frac{\pi ^4}{18 x^3}+\frac{\pi ^6}{54 x^5}-\frac{5 \pi ^8}{648 x^7}+\frac{7 \pi ^{10}}{1944 x^9}-\frac{7 \pi ^{12}}{3888 x^{11}}+O\left(\frac{1}{x^{12}}\right),\\[9pt]
    &z=x+\frac{2 \pi ^2}{3 x}-\frac{2 \pi ^4}{9 x^3}+\frac{4 \pi ^6}{27 x^5}-\frac{10 \pi ^8}{81 x^7}+\frac{28 \pi ^{10}}{243 x^9}-\frac{28 \pi ^{12}}{243 x^{11}}+O\left(\frac{1}{x^{12}}\right),\\[9pt]
    &w=x+\frac{\pi ^2}{x}-\frac{\pi ^4}{2 x^3}+\frac{\pi ^6}{2 x^5}-\frac{5 \pi ^8}{8 x^7}+\frac{7 \pi ^{10}}{8 x^9}-\frac{21 \pi ^{12}}{16 x^{11}}+O\left(\frac{1}{x^{12}}\right).
\end{align*}
It is readily checked that for $x\geq 4$,
\begin{align}
    &y_1<y<y_2,\label{Ineq-y}\\[9pt]
    &z_1<z<z_2,\label{Ineq-z}\\[9pt]
    &w_1<w<w_2,\label{Ineq-w}
\end{align}
where
\begin{align*}
    &y_1=x+\frac{\pi ^2}{3 x}-\frac{\pi ^4}{18 x^3}+\frac{\pi ^6}{54 x^5}-\frac{5 \pi ^8}{648 x^7}+\frac{7 \pi ^{10}}{1944 x^9}-\frac{7 \pi ^{12}}{3888 x^{11}},\\[9pt]
    &y_2=x+\frac{\pi ^2}{3 x}-\frac{\pi ^4}{18 x^3}+\frac{\pi ^6}{54 x^5}-\frac{5 \pi ^8}{648 x^7}+\frac{7 \pi ^{10}}{1944 x^9},\\[9pt]
    &z_1=x+\frac{2 \pi ^2}{3 x}-\frac{2 \pi ^4}{9 x^3}+\frac{4 \pi ^6}{27 x^5}-\frac{10 \pi ^8}{81 x^7}+\frac{28 \pi ^{10}}{243 x^9}-\frac{28 \pi ^{12}}{243 x^{11}},\\[9pt]
    &z_2=x+\frac{2 \pi ^2}{3 x}-\frac{2 \pi ^4}{9 x^3}+\frac{4 \pi ^6}{27 x^5}-\frac{10 \pi ^8}{81 x^7}+\frac{28 \pi ^{10}}{243 x^9},\\[9pt]
    &w_1=x+\frac{\pi ^2}{x}-\frac{\pi ^4}{2 x^3}+\frac{\pi ^6}{2 x^5}-\frac{5 \pi ^8}{8 x^7}+\frac{7 \pi ^{10}}{8 x^9}-\frac{21 \pi ^{12}}{16 x^{11}},\\[9pt]
    &w_2=x+\frac{\pi ^2}{x}-\frac{\pi ^4}{2 x^3}+\frac{\pi ^6}{2 x^5}-\frac{5 \pi ^8}{8 x^7}+\frac{7 \pi ^{10}}{8 x^9}.
\end{align*}

With these bounds of $y,z$ and $w$ in \eqref{Ineq-y}, \eqref{Ineq-z} and \eqref{Ineq-w}, we are now in a position to estimate $t_1, t_2, t_3, e^{w+y-2z}$ and $e^{x-2y+z}$ in terms of $x$.

First, we consider $t_1, t_2$, and $t_3$.
By the definition of $\alpha(t)$,
\begin{align*}
     \alpha(w)=w^{10}-w^{9}+1.
\end{align*}
Noting that $w^9=(x^2+2\pi^2)^4\sqrt{x^2+2\pi^2}$, which involves a radical, to give a
feasible estimate for $w^9$ without a radical, we may make use of \eqref{Ineq-w} to deduce that for $x\geq 4$,
\begin{equation*}
    w_1 w^8<w^9<w_2 w^8.
\end{equation*}
Let
\begin{equation*}
    \eta_1=w^{10}-w_1 w^{8}+1,
\end{equation*}
so that for $x\geq 4$,
\begin{equation}\label{eta-1}
     \alpha(w)<\eta_1.
\end{equation}
Similarly, set
\begin{align*}
\eta_2=&y^{30}-3 y_1 y^{28}+3 y^{28}-y_1 y^{26}+3 y^{20}-6 y_1 y^{18}+3 y^{18}+3 y^{10}-3 y_1 y^8+1,\\[9pt]
\eta_3=&z^{30}-3 z_2 z^{28}+3 z^{28}-z_2 z^{26}-3 z^{20}+6 z_1 z^{18}-3 z^{18}+3 z^{10}-3 z_2 z^8-1,\\[9pt]
\eta_4=&y^{20}-2 y_2 y^{18}+y^{18}+2 y^{10}-2 y_2 y^8+1,\\[9pt]
\eta_5=&z^{20}-2 z_2 z^{18}+z^{18}-2 z^{10}+2 z_1 z^8+1.
\end{align*}
Then we have for $x\geq 4$,
\begin{equation}\label{eta2-5}
    \alpha(y)^3<\eta_2,\ \beta(z)^3>\eta_3,\ \alpha(y)^2>\eta_4,\ \beta(z)^2>\eta_5.
\end{equation}

Employing the relations in \eqref{eta-1} and \eqref{eta2-5}, we deduce that for $x\geq 4$,
\begin{align}
&t_1= x^{12} z^{36} \alpha(y)^3\alpha(w)< x^{12} z^{36} \eta_1 \eta_2,\label{t-1-1}\\[9pt]
&t_2= (x^{10}-x^9-1)y^{36} w^{12} \beta(z)^3> (x^{10}-x^9-1)y^{36}w^{12} \eta_3,\label{t-2-1}\\[9pt]
&t_3= x^7 y^{12} z^{12} w^{12} \alpha(y)^2\beta(z)^2> x^7 y^{12} z^{12} w^{12}\eta_4 \eta_5.\label{t-3-1}
\end{align}

We continue to estimate $e^{w+y-2z}$ and $e^{z+ x-2y}$.
Applying \eqref{Ineq-y}, \eqref{Ineq-z} and \eqref{Ineq-w} to $w+y-2z$, we see that for $x\geq 4$,
\begin{align}
w+y-2z<w_2+y_2-2z_1,
\end{align}
which implies that
\begin{equation}\label{ineq-E-7}
    e^{w+y-2z}<e^{w_2+ y_2-2z_1}.
\end{equation}
In order to give a feasible upper bound for $e^{w+y-2z}$, we define
\begin{equation}\label{M}
    \Phi(t)=1+t +\frac{t ^2}{2}+\frac{t ^3}{6}+\frac{t ^4}{24}+\frac{t ^5}{120}+\frac{t^6}{720},
\end{equation}
and it can be proved that for $t<0$,
\begin{equation}\label{in-E-M}
 e^{t}<\Phi(t).
\end{equation}
Note that
\begin{equation*}
w_2+y_2-2z_1=-\frac{\pi ^4 (108 x^8-216 \pi ^2 x^6+375 \pi ^4 x^4-630 \pi ^6 x^2-224 \pi ^8)}{972 x^{11}}<0
\end{equation*}
holds for $x\geq 5$, since
\begin{equation*}
108 x^8-216 \pi ^2 x^6>0
\end{equation*}
for $ x>\sqrt{2} \pi \approx4.443 $, and
\begin{equation*}
375 \pi ^4 x^4-630 \pi ^6 x^2-224 \pi ^8>0
\end{equation*}
for $x> \frac{\pi}{5}\sqrt{\sqrt{{2443}/{3}}+21}\approx 4.422$.
Thus, by \eqref{in-E-M}, we obtain that for $x\geq 5$,
\begin{eqnarray}\label{Ineq-E-8}
e^{w_2+ y_2-2z_1}<\Phi(w_2+y_2-2z_1).
\end{eqnarray}
Combining \eqref{ineq-E-7} and \eqref{Ineq-E-8} yields that for $x\geq 5$,
\begin{eqnarray}\label{e1}
e^{w+ y-2z}<\Phi(w_2+y_2-2z_1).
\end{eqnarray}

Similarly, applying \eqref{Ineq-y}, \eqref{Ineq-z} and \eqref{Ineq-w} to $z+x-2y$, we find that for $x\geq 4$,
\begin{equation}\label{E-2}
   z_1+x-2 y_2< z+x-2y,
\end{equation}
so that
\begin{equation}\label{Ineq-E-9}
e^{z+x-2y}>e^{z_1+x-2y_2}.
\end{equation}
Define
\begin{equation}\label{m}
 \phi(t)=1+t+\frac{t ^2}{2}+\frac{t ^3}{6}+\frac{t ^4}{24}+\frac{t ^5}{120}+\frac{t^6}{720}+\frac{t^7}{5040}.
\end{equation}
It is true that for $t<0$,
\begin{equation}\label{in-E-m}
 \phi(t)<e^{t}.
\end{equation}
We now give a lower bound for $e^{z_1+x-2y_2}$. Since 
\begin{align*}
 z+x-2y&=\sqrt{x^2+\frac{4\pi^2}{3}}+x-2\sqrt{x^2+\frac{2\pi^2}{3}}\\[9pt]
 &=\frac{-\left(\sqrt{x^2+4\pi^2/3}-x\right)^2}
 {\sqrt{x^2+{4\pi^2}/{3}}+x+2\sqrt{x^2+{2\pi^2}/{3}}},
\end{align*}
which is negative for $n\geq 2$, by \eqref{E-2}, we deduce that for  $x\geq 4$,
\begin{equation}\label{MMM}
    z_1+x-2y_2<0.
\end{equation}
Thus,  applying  \eqref{in-E-m} to \eqref{MMM} gives us that  for $x\geq4$,
\begin{eqnarray}\label{Ineq-E-10}
e^{z_1+x-2y_2}>\phi(z_1+x-2y_2).
\end{eqnarray}
Combining \eqref{Ineq-E-9} and \eqref{Ineq-E-10}, we find that for $x\geq4$,
\begin{eqnarray}\label{e2}
e^{z+x-2y}>\phi(z_1+x-2y_2).
\end{eqnarray}

Using the above bounds for  $t_1, t_2, t_3, e^{w+z-2y}$ and $e^{z+x-2y}$,
we obtain that for $x\geq5$,
\begin{align}
    & -e^{w+y-2z}t_1+e^{z+ x-2y}t_2+110t_3\nonumber \\[9pt]
    &\quad >-\Phi(w_2+y_2-z_1)x^{12} z^{36} \eta_1 \eta_2
        +\phi(z_1+x-2y_2)(x^{10}-x^9-1)y^{36}w^{12} \eta_3\nonumber \\[9pt]
    &\qquad \quad +110 x^7 y^{12} z^{12} w^{12}\eta_4 \eta_5.\label{ineq4-2}
\end{align}
To verify \eqref{eq4.1}, we show that for $x\geq 358$,
\begin{align}\label{ineq-END}
  &-\Phi(w_2+y_2-z_1)x^{12} z^{36} \eta_1 \eta_2+\phi(z_1+x-2y_2)(x^{10}-x^9-1)y^{36}w^{12} \eta_3\nonumber \\[9pt]
  &\qquad+110 x^7 y^{12} z^{12} w^{12}\eta_4 \eta_5>0.
\end{align}
 Substituting $y,z$ and $w$ with $\sqrt{x^2+2\pi^2/3},\sqrt{x^2+4\pi^2/3}$ and $\sqrt{x^2+2\pi^2}$ respectively,
the left hand side of the inequality \eqref{ineq-END} can be expressed as $H(x)/ G(x)$, where
\[H(x)=\sum_{k=0}^{171} a_k x^k\]
    and
\[
    G(x)=39686201656473354776757087428535162639482880 x^{88}.
\]
Here we just list the values of $a_{169}, a_{170}$ and $a_{171}$:
\begin{align*}
    &a_{169}=734929660305062125495501619046947456286720 \\[9pt]
    &\qquad \qquad \times\left(35640+261360 \pi ^2-194 \pi ^6-249 \pi ^8\right),\\[9pt]
    &a_{170}=5879437282440497003964012952375579650293760 \left(7 \pi ^6-2970\right),\\[9pt]
    &a_{171}=4409577961830372752973009714281684737720320 \left(990-\pi ^6\right),
\end{align*}
which are all positive.

Given that $G(x)$ is always positive, we aim to prove that $H(x)>0$.
Apparently, $x\geq 2$ for $n\geq 2$ and hence
\begin{equation}\label{ineqH1}
  H(x) \geq \sum_{k=0}^{170} -|a_k| x^k+a_{171}x^{171}.
\end{equation}
Moreover, numerical evidence indicates that for  any $0\leq k\leq 168$,
\begin{equation}\label{ineq-XX}
   -|a_k| x^k>-a_{169} x^{169}
\end{equation}
holds for $x\geq 181$. It follows that
 for $x\geq 181$,
\begin{equation}\label{In-p-1}
    \sum_{k=0}^{170} -|a_k| x^k+a_{171}x^{171}>
    (-170 a_{169}-a_{170}x+a_{171}x^2)x^{169}.
\end{equation}
Combining \eqref{ineqH1} and \eqref{In-p-1}, we obtain that for $x\geq 181$,
\begin{equation}
    H(x)> (-170 a_{169}-a_{170} x+a_{171}x^2)x^{169}.
\end{equation}
Thus, $H(x)$ is positive provided
\begin{equation}
    -170 a_{169}-a_{170} x+a_{171}x^2>0,
\end{equation}
which is true if
\begin{equation*}
    x>\frac{\sqrt{a_{170}^2+680 a_{169} a_{171}}+a_{170}}{2 a_{171}}\approx  357.867.
\end{equation*}
Hence we conclude that $H(x)$ is positive when $x\geq 358$. This
 proves \eqref{ineq-END}.

Combining \eqref{ineq4-2} and \eqref{ineq-END}, we find that for $x\geq 358$, or equivalently, for $n\geq 19480$, \eqref{eq4.1} holds, that is,
\begin{equation}\label{xii}
 -e^{w+y-2z}t_1+e^{z+ x-2y}t_2+110t_3>0.
 \end{equation}
For $2\leq n\leq 19480$,  \eqref{xii} can be directly verified.
This completes the proof.
\qed

\section{An inequality on $u_n$ and $f(n)$}

In this section, we present an inequality on
$u_n$ and $f(n)$ that is also needed in the proof of
Theorem \ref{Theorem 2.2}.

\begin{theo} \label{Lemmax5}
Let $u_n$ be defined as \eqref{1.8}. For
 $0<t<1$, let
\begin{equation}\label{Q}
 Q(t)=\frac{3 t +2 \sqrt{(1-t)^3}-2}{t^2}.
 \end{equation}
Then for $n\geq 85$,
\begin{align}\label{Ineq-The-4-1-XXX}
f(n)+\frac{110}{\mu(n-1)^5}< Q(u_n).
\end{align}
\end{theo}

The proof of this theorem is based on the following Lemma, which gives an upper bound of $f(n)$.
Recall that
\[f(n)=e^{x-2 y+z}\frac{\left(x^{10}-x^9-1\right) y^{24} \left(z^{10}-z^9-1\right)}{x^{12} \left(y^{10}-y^9+1\right)^2 z^{12}},\]
where $x,y,z,w$ are defined as in \eqref{xyzw}.

\begin{lemm}\label{Lemma-4-2}
Let
\[\Phi(t)=1+t +\frac{t ^2}{2}+\frac{t ^3}{6}+\frac{t ^4}{24}+\frac{t ^5}{120}+\frac{t^6}{720},\]
as defined in \eqref{M}, and let $y_1, y_2, z_1$ and $z_2$ be defined as in the proof of Theorem \ref{Lemmax2}.
For $n\geq 4$, we have
\begin{equation}\label{Ineq-f}
f(n)<\frac{\Phi(x-2y_1+z_2) y^{24} (x^{10}-x^9-1)\left(z^{10}-z^{8} z_1-1\right)}{x^{12} z^{12}(y^{20}-2 y^{18} y_2+y^{18}+2 y^{10}-2 y^8y_2+1)}<1.
\end{equation}
\end{lemm}

\proof Using the bounds for $y,z, w$ as given in \eqref{Ineq-y}, \eqref{Ineq-z} and \eqref{Ineq-w}, we  shall derive  estimates for the factors
 $e^{x-2 y+z}$,
$z^{10}-z^9-1$, and $\left(y^{10}-y^9+1\right)^2$ that appear in $f(n)$.
It is easily verified that for $x\geq  4$,
 \begin{gather}
    e^{x-2y+z}<e^{x-2y_1+z_2},\label{InE-4}\\[9pt]
    z^{10}-z^9-1<z^{10}-z^{8} z_1-1,\label{Ineq-z-4}\\[9pt]
    \left(y^{10}-y^9+1\right)^2>y^{20}-2 y^{18} y_2+y^{18}+2 y^{10}-2 y^8y_2+1. \label{Ineq-y-4}
 \end{gather}

We further give an   upper bound for $e^{x-2y_1+z_2}$.
Write
\begin{equation}\label{x-y-z}
    x-2y_1+z_2=-\frac{\pi ^4 \left(216 x^8-216 \pi ^2 x^6+210 \pi ^4 x^4-210 \pi ^6 x^2-7 \pi ^8\right)}{1944 x^{11}}.
\end{equation}
For $ x>\pi$, we have
\begin{equation*}
216 x^8-216 \pi ^2 x^6>0,
\end{equation*}
and for $x>\pi\sqrt{(\sqrt{17/15}+1)/2}\approx 3.192$, we have
\begin{equation*}
210 \pi ^4 x^4-210 \pi ^6 x^2-7 \pi ^8>0.
\end{equation*}
Therefore, it follows from \eqref{x-y-z}   that for $x\geq 4$,
\[x-2y_1+z_2<0,\]
which, together with \eqref{in-E-M}, yields that for $x\geq 4$,
\begin{equation}\label{InE-4-2}
e^{x-2y_1+z_2}<\Phi(x-2y_1+z_2).
\end{equation}
Combining \eqref{InE-4} and \eqref{InE-4-2}, we find that for $x\geq 4$,
\begin{equation}\label{h2}
e^{x-2y+z}<\Phi(x-2y_1+z_2).
\end{equation}
By means of the estimates in  \eqref{Ineq-z-4}, \eqref{Ineq-y-4} and \eqref{h2},
we arrive at the first inequality in  \eqref{Ineq-f}.

To prove the second inequality in \eqref{Ineq-f}, by expressing   $y,z$ and $w$
in terms of $x$, we get
\begin{equation*}
    \frac{\Phi(x-2y_1+z_2) y^{24} (x^{10}-x^9-1)\left(z^{10}-z^{8} z_1-1\right)}{x^{12} z^{12}(y^{20}-2 y^{18} y_2+y^{18}+2 y^{10}-2 y^8y_2+1)}=\frac{H(x)}{G(x)},
\end{equation*}
 where $H(x)$ and $G(x)$ are both polynomials of degree $121$. Write
\begin{equation}\label{Eq-H-G}
    H(x)=\sum_{k=0}^{121} b_k x^k,\  G(x)=\sum_{k=0}^{121} c_k x^{k}.
\end{equation}
Here are the values of $b_k$ and $c_k$ for $116\leq k\leq 121$:
\begin{align*}
&b_{116}=-1398983398232765780459520 \pi ^4 \left(5181+41 \pi ^2\right),\\[9pt]
&b_{117}=25181701168189784048271360 \pi ^2 \left(21+151 \pi ^2\right),\\[9pt]
&b_{118}=-4196950194698297341378560 \pi ^2 \left(258+\pi ^2\right),\\[9pt]
&c_{116}=-7197769583907579940464230400 \pi ^4,\\[9pt]
&c_{117}=75545103504569352144814080 \pi ^2 \left(7+50 \pi ^2\right),\\[9pt]
&c_{118}=-1082813150232160714075668480 \pi ^2,\\[9pt]
    &b_{119}=c_{119}=12590850584094892024135680 \left(3+44 \pi ^2\right),\\[9pt]
    &b_{120}=c_{120}=-75545103504569352144814080,\\[9pt]
    &b_{121}=c_{121}=37772551752284676072407040.
\end{align*}

For our purpose, we claim that for $x\geq 135$,
\begin{equation}\label{In-G}
    G(x)>0,
\end{equation}
and
\begin{equation}\label{H-G}
    G(x)-H(x)>0.
\end{equation}

Observe that for  $0\leq k\leq 118$,
\begin{equation}
-|c_k| x^k>-c_{119} x^{119}
\end{equation}
holds when \[ x> \pi  \sqrt{\frac{6 \left(7+50 \pi ^2\right)}{3+44 \pi ^2}}\approx8.232.\]
It follows that for $x\geq 9$,
\begin{equation}\label{G-2}
    G(x)>(-120 c_{119}+ c_{120} x +c_{121} x^2)x^{119}.
\end{equation}
Since
\begin{equation}
    120 c_{119}+ c_{120} x +c_{121} x^2>0
\end{equation}
  for \[ x>1+\sqrt{11 \left(11+160 \pi ^2\right)}\approx 133.255,\] we
obtain that $G(x)>0$ for $x\geq 134$.

Similarly, to prove \eqref{H-G}, we observe that for $0\leq k\leq 115$,
\begin{equation}
-|c_k-b_k| x^k>-(c_{116}-b_{116}) x^{116}
\end{equation}
 for \[ x>\frac{1}{2} \pi  \sqrt{\frac{5616+3127 \pi ^2}{108+123 \pi ^2}}\approx 8.232.\]
Hence
\begin{align}
    &G(x)-H(x)=\sum_{k=0}^{118}(c_k-b_k) x^{118}\nonumber\\[5pt]
    &\qquad >(-117(c_{116}-b_{116})+(c_{117}-b_{117}) x + (c_{118}-b_{118}) x^2)x^{118}\label{H-G-2},
\end{align}
which is positive for $n\geq 135$, here we have used the fact that
\begin{equation}
    -117(c_{116}-b_{116})+(c_{117}-b_{117}) x + (c_{118}-b_{118}) x^2>0
\end{equation}
if  \[ x>3+\sqrt{3 \left(471+533 \pi ^2\right)}\approx 134.128.\]

Combining \eqref{In-G} and \eqref{H-G}, we deduce that the second inequality \eqref{Ineq-f} is valid for $x\geq 135$, or equivalently, for $n\geq 2771$.
The case for $4\leq n \leq 2771$ can be directly verified, and hence the proof
is complete.
\qed

We are now ready to prove Theorem \ref{Lemmax5}.

\noindent
{\it Proof of Theorem \ref{Lemmax5}.} Recall that the theorem
states that for $n\geq 85$, 
\begin{align}\label{Ineq-The-4-1}
f(n)+\frac{110}{\mu(n-1)^5}< Q(u_n).
\end{align}
It can be checked that  \eqref{Ineq-The-4-1} is true   for $85 \leq n \leq 35456$.
We now show that \eqref{Ineq-The-4-1} is true for $n\geq 35457$.
Since
\begin{equation}
Q'(t)=\frac{1}{\left(\sqrt{1-t }+1\right)^3},
\end{equation}
which is positive  for $0<t<1$,
 $Q(t)$ is  increasing for $0<t<1$.
 By Theorem \ref{Lemmax1}, we know that  $f(n)<u_n$ for $n\geq 1207$,
 so that for $n\geq 1207$,
 \begin{equation}\label{Q(u)}
    Q(f(n))<Q(u_n).
 \end{equation}
Thus \eqref{Ineq-The-4-1} is justified if we can prove that for $n\geq 35457$,
\begin{align}\label{2.3}
f(n)+\frac{110}{\mu(n-1)^5}< Q(f(n)).
	\end{align}

Let
\begin{align}\label{psi-t}
\psi(t)=Q(t)-t=\frac{3 t+2 \sqrt{(1-t )^3}-t ^3-2}{t^2}.
\end{align}
In this notation, \eqref{2.3} says that for $n\geq 35457$,
\begin{eqnarray}\label{fk1}
\psi(f(n))>\frac{110}{\mu(n-1)^5}.
\end{eqnarray}

To prove the above inequality, we shall use the polynomials
$G(x)$ and $H(x)$ as given by \eqref{Eq-H-G}. More specifically,
\[H(x)=\sum_{k=0}^{121} b_k x^k,\  G(x)=\sum_{k=0}^{121} c_k x^{k}.\]
Note that $\psi(t)$ is decreasing for $0<t<1$, since for $0<t<1$,
\begin{eqnarray*}
\psi'(t)=-\frac{\sqrt{1-t } \left(-t +3 \sqrt{1-t }+4\right)}{\left(\sqrt{1-t }+1\right)^3}<0.
\end{eqnarray*}

It can be seen from Lemma \ref{Lemma-4-2}  that $0<f(n)<{H(x)}/{G(x)}<1$ for $n\geq 4$, so that for $n \geq 35457$,
\begin{equation}\label{Ineq-H-4}
    \psi(f(n))>\psi\left(\frac{H(x)}{G(x)}\right).
\end{equation}
 Because of (\ref{Ineq-H-4}), to verify (\ref{fk1}), it is sufficient to show that for $n\geq 35457$,
\begin{equation}\label{Ineq-H-3-G}
    \psi\left(\frac{H(x)}{G(x)}\right)
  >\frac{110}{\mu(n-1)^{5}}.
\end{equation}
This goal can be achieved by
 finding an estimate for $\psi\left({H(x)}/{G(x)}\right)$.
 We first derive the following range of ${H(x)}/{G(x)}$ for $x\geq 134$,
\begin{equation}\label{h11}
\frac{\sqrt{5}-1}{2}<\frac{H(x)}{G(x)}<1.
\end{equation}
By Lemma \eqref{Lemma-4-2}, we know that ${H(x)}/{G(x)}<1$ for $x \geq 4$ and $G(x)>0$ for $x\geq 134$. To justify \eqref{h11}, we  only need to show that for $x\geq 134$,
\begin{equation}\label{H2}
    2H(x)-(\sqrt{5}-1)G(x)>0.
\end{equation}
Note that
\[b_{119}=c_{119},\qquad b_{120}=c_{120}, \qquad b_{121}=c_{121},\]
and observe that for $0\leq k \leq 118$,
\begin{equation}
    -|2b_k-(\sqrt{5}-1)c_k| x^{k}>-(3 -\sqrt{5})c_{119} x^{119} 
\end{equation}
when
\[x>\pi  \sqrt{\frac{\pi ^2 \left(\sqrt{5}+303\right)+42}{3+44 \pi ^2}}\approx8.303.\]
It follows that for $x\geq 9$,
\begin{equation}
    2H(x)-(\sqrt{5}-1)G(x)>(3 -\sqrt{5})(-120c_{119}+c_{120} x +c_{121} x^2)x^{119}.
\end{equation}
Since
\begin{equation*}
    -120c_{119}+c_{120} x +c_{121} x^2>0
\end{equation*}
for $ x>\sqrt{11 \left(11+160 \pi ^2\right)}+1\approx 133.255$,
we arrive at \eqref{H2}, and so  \eqref{h11} is proved.

The above range of $H(x)/G(x)$ enables us to bound $\psi\left({H(x)}/{G(x)}\right)$.
It turns out that in the same range $\frac{\sqrt{5}-1}{2} <t <1$, we have \begin{eqnarray}\label{psi}
\psi(t)>(1-t)^{3/2}, 	
\end{eqnarray}
since
\[\psi(t)-(1-t)^{3/2}=\frac{(1-t )^{3/2} \left(t+\frac{\sqrt{5}+1}{2} \right)\left(t-\frac{\sqrt{5}-1}{2} \right)}{\left(\sqrt{1-t }+1\right)^2(\sqrt{1-t}+t)},\]
which is positive for $\frac{\sqrt{5}-1}{2} <t <1$.

In view of  \eqref{h11} and \eqref{psi}, we infer that for $x\geq 134$,
\begin{equation}\label{psi-H-G}
    \psi\left(\frac{H(x)}{G(x)}\right)
    >\left(1-\frac{H(x)}{G(x)}\right)^{\frac{3}{2}}.
\end{equation}
 We continue to  show that for $x\geq 483$, or equivalently, for $n\geq 35457$,
\begin{equation}\label{Ineq-H-4-G}
   \left(1-\frac{H(x)}{G(x)}\right)^{\frac{3}{2}}>\frac{110}{\mu(n-1)^{5}}.
\end{equation}
Since $G(x)>0$ for $x\geq 134$, the above inequality can be reformulated
 as follows. For $x\geq 483$,
\begin{equation}\label{H3}
x^{10}(G(x)-H(x))^3 -110^2 G(x)^3 >0.
\end{equation}
The left hand side of (\ref{H3}) is a polynomial of degree 364, and 
we write
\begin{equation}
x^{10}(G(x)-H(x))^3 -110^2 G(x)^3=\sum_{k=0}^{364} \gamma_k x^{k}.
\end{equation}
The values of $\gamma_{364}, \gamma_{363}$ and $\gamma_{362}$ are given below:
\begin{align*}
    &\gamma_{364}=2^{72} 3^{105} 5^3 \pi^{12},\\[9pt]
    &\gamma_{363}=-2^{73} 3^{107} 5^3 \left(490050+\pi ^{12}\right),\\[9pt]
    &\gamma_{362}=2^{72} 3^{105} 5^3 \left(52925400+144 \pi ^{12}+41 \pi ^{14}\right).
\end{align*}
For $0\leq k \leq 361$, we find that
\begin{equation}
    -|\gamma_k|x^k>-\gamma_{362} x^{362},
\end{equation}
provided that
\[x>\frac{793881000+2328717600 \pi ^2+3996 \pi ^{12}+4392 \pi ^{14}+\pi ^{16}}{317552400+864 \pi ^{12}+246 \pi ^{14}}\approx20.126.\]
Thus, for $x\geq 21$,
\begin{align}
x^{10}(G(x)-H(x))^3 -110^2 G(x)^3>(-363 \gamma_{362}+\gamma_{363} x + \gamma_{364} x^2)x^{362},
\end{align}
which is positive, since
\begin{equation*}
    -363 \gamma_{362}+\gamma_{363} x + \gamma_{364} x^2>0
\end{equation*}
as long as
\[x>\frac{\sqrt{1452 \gamma_{362} \gamma_{364}+\gamma_{363}^2}-\gamma_{363}}{2 \gamma_{364}}\approx 482.959.\]
Hence \eqref{Ineq-H-4-G} is confirmed. Combining \eqref{psi-H-G} and \eqref{Ineq-H-4-G}, we are led to \eqref{Ineq-H-3-G}. The proof is therefore complete.
\qed

\section{Proof of Theorem \ref{Theorem 2.2}}

In this section, we present a proof of   Theorem \ref{Theorem 2.2}
based on the intermediate inequalities in the previous sections.
 The theorem states that for $n\geq 95$.
\begin{equation}\label{U-n}
    4(1-u_n)(1-u_{n+1})-(1-u_n u_{n+1})^2>0,
\end{equation} where
\[u_n=\frac{p(n+1)p(n-1)}{p(n)^2}.\]

\noindent
{\it Proof of Theorem \ref{Theorem 2.2}.}
We shall make use of the fact that $u_n<1$ for $n\geq 26$, as proved by
DeSalvo and Pak \cite{Desalvo}. In order to prove \eqref{U-n},
we define $F(t)$ to be a function in $t$:
\begin{align}\label{Fun-F}
    F(t)&=4 (1-u_n)(1-t)-(1- u_n t)^2.
\end{align}
Then \eqref{U-n} says that for $n\geq 95$,
\begin{equation}\label{In-5-F}
    F(u_{n+1})>0.
\end{equation}
For $95 \leq n \leq 1206$, \eqref{In-5-F} can be directly checked.
We proceed to prove that \eqref{In-5-F}  holds for $n\geq 1207$.
Let $Q(t)$ be as defined in \eqref{Q}, that is,
\[Q(t)=\frac{3 t +2 \sqrt{(1-t)^3}-2}{t^2}.\]
We claim that  $F(t)>0$ for $u_n<t<Q(u_n)$. Rewrite $F(t)$ as
\[F(t)=-u_n^2 t^2+(6 u_n-4) t-4 u_n+3.\]
The equation $F(t)=0$ has two solutions:
\begin{equation*}
    P(u_n)=\frac{3 u_n -2 \sqrt{(1-u_n)^3}-2}{u_n ^2},\quad \ Q(u_n)=\frac{3 u_n +2 \sqrt{(1-u_n)^3}-2}{u_n^2},
\end{equation*}
so that $F(t)>0$ for $P(u_n)<t<Q(u_n)$.
Furthermore, we see that
\[F(u_n)=(1-u_n)^3 (u_n +3)>0,\]
which implies $P(u_n)<u_n<Q(u_n)$.
Therefore, $F(t)>0$ for $u_n<t<Q(u_n)$, as claimed.

To obtain \eqref{In-5-F}, it remains to show that for $ n\geq 1207 $,
\begin{equation}\label{LastIneq}
    u_n<u_{n+1}<Q(u_n).
\end{equation}
Recall that $u_n < u_{n+1}$  holds for $n\geq 116$, as
proved  by Chen, Wang and Xie \cite{Chen2}.
By Theorem \ref{Lemmax1}, we know that $u_{n+1}<g(n+1)$ for $n\geq 1207$.
 But Theorem \ref{Lemmax2} asserts that for $n\geq 2$,
\[g(n+1)<f(n)+\frac{110}{\mu(n-1)^5}.\] 
Furthermore, Theorem \ref{Lemmax5} states that for $n\geq 2$,
\[f(n)+\frac{110}{\mu(n-1)^5}<Q(u_n).\]
Thus we conclude that $u_{n+1}<Q(u_n)$ for $n\geq 1207$, as required.
\qed

\end{document}